%% file: main.tex
\documentclass[runningheads]{llncs}
\usepackage[T1]{fontenc}
\usepackage{graphicx}
\usepackage{amsmath}
\usepackage{amssymb}
\usepackage{booktabs}
\usepackage{siunitx}
\usepackage[table]{xcolor}
\definecolor{matbetter1}{RGB}{220, 245, 220}   
\definecolor{matbetter2}{RGB}{190, 235, 190}   
\definecolor{matbetter3}{RGB}{160, 225, 160}   
\definecolor{matworse1}{RGB}{255, 220, 220}    
\definecolor{matworse2}{RGB}{250, 195, 195}    
\definecolor{matworse3}{RGB}{245, 170, 170}    
\definecolor{matequiv1}{RGB}{220, 235, 250}    
\definecolor{matequiv2}{RGB}{195, 220, 245}    
\definecolor{matequiv3}{RGB}{170, 205, 240}    
\usepackage[
    colorlinks=true,
    linkcolor=blue,
    citecolor=blue,
    urlcolor=blue
]{hyperref}
\begin{document}
\title{Impact of Scaling and Rounding on Metaheuristic Performance for the Vehicle Routing Problem with Time Windows}

\titlerunning{Scaling and Rounding VRPTW Benchmarks}
%
\author{Florian Rascoussier\inst{1,2}\orcidID{0009-0005-3253-9814} \and
Romain Billot\inst{1} \and
Lina Fahed\inst{1}  \and
Christine Solnon\inst{2}}
\authorrunning{F. Rascoussier et al.}
%
\institute{IMT Atlantique, Technop\^ole Brest-Iroise, 29238 Brest, France \and
CITI, INSA Lyon, INRIA, France\\
\email{florian.rascoussier@insa-lyon.fr}}
\maketitle              
\begin{abstract}
Classical Euclidean instances for the Vehicle Routing Problem with Time Windows (VRPTW) have floating-point arc costs derived from node coordinates. This raises reproducibility and numerical consistency issues. Hence, a common practice is to scale and round data to integer values, yet the impact of these choices is poorly understood. This paper studies the effect of integer scaling by a factor $P$ and compares three rounding schemes with distinct feasibility and optimality guarantees. We analyze how scaling and rounding influence empirical performance of two well-known VRPTW solvers: Hybrid Genetic Search (HGS) and OR-Tools.
We evaluate the quality of solutions mapped back to the original floating-point instances. Our results show that scaling and rounding can significantly affect both solver performance and solution quality. Based on these findings, we provide practical recommendations for scaling and rounding in VRPTW benchmarks to improve efficiency, robustness and reproducibility.

\end{abstract}
\section{Introduction}
\label{sec:intro}

The Vehicle Routing Problem with Time Windows (VRPTW) is a cornerstone of combinatorial optimization. The most widely used benchmark instances, such as those by Solomon \cite{solomonAlgorithmsVehicleRouting1987} and Gehring \& Homberger \cite{HombergerGehring1999}, are defined on a 2D plane where arc costs correspond to Euclidean distances. These distances are inherently floating-point numbers as an approximation of reals. However, solvers often require or prefer integer values for performance, reproducibility, and numerical stability.

Translating floating-point values into integers involves two key operations: \emph{scaling} the values by a factor $P$ (typically a power of 10) to preserve precision and \emph{rounding} the scaled values to near integers. While these steps are often treated as minor preprocessing details, they fundamentally alter instances: 
different choices of $P$ and rounding strategies can transform the feasible region and fitness landscape, potentially impacting solver performance and optimal solutions. For instance, the 2021 DIMACS challenge \cite{DIMACS2021} adopted a scaled and truncated approach, whereas the classic SINTEF Best-Known Solutions (BKS) \cite{Sintef2008} are evaluated on the original floating-point instances; a divergence that complicates direct comparisons and can alter conclusions about solver performance.

In this paper, we systematically investigate this impact using two well-known open-source solvers: Hybrid Genetic Search (HGS) \cite{vidalHybridGeneticSearch2022} and OR-Tools \cite{&ortools}. 
Our contributions are:
(i)~a formal characterization of three rounding schemes as inner and outer approximations of VRPTW instances, with distinct feasibility and optimality guarantees;
(ii)~a large-scale empirical analysis quantifying their impact on feasibility transfer, objective distortion, and convergence behavior;
(iii)~a statistical demonstration of evaluation-dependent ranking inversion between preprocessing configurations;
(iv)~practical recommendations for reproducible and comparable VRPTW benchmarking.

The remainder of this paper is organized as follows. In Section~\ref{sec:back}, we define the VRPTW and introduce notations. In Section~\ref{sec:methodology}, we present the scaling and rounding schemes. Our experimental setup is detailed in Section~\ref{sec:experiments}. Results are discussed in Section~\ref{sec:results}. In Section~\ref{sec:conclusion}, we conclude and outline perspectives.

\section{Background and Notations}\label{sec:back}

A Euclidean VRPTW instance is defined by a vehicle capacity $\kappa$, a set $C$ of $n$ customer points to visit and a depot $d$. For each point $i \in C\cup\{d\}$, $(x_i, y_i)$ denotes its coordinates, and for each customer point $i\in C$, $q_i$ denotes its demand, $[a_i, b_i]$ its time window (TW), and $s_i$ its service time. The travel time $\tau_{ij}$ between two points $i$ and $j$ is given by their Euclidean distance, {\em i.e.} $\tau_{ij} = \sqrt{(x_i - x_j)^2 + (y_i - y_j)^2}$. A feasible solution $S$ is a set of $r$ tours $S = \{T_1, \ldots, T_r\}$ such that (i) every point in $C$ belongs to exactly one tour, (ii) each tour $T_v \in S$ starts and ends at the depot, (iii) vehicle capacity constraints are satisfied, {\em i.e.}, the sum of demands of the customer points does not exceed $\kappa$, and (iv) time window constraints are satisfied, {\em i.e.}, each point $i$ is serviced within its time window $[a_i, b_i]$ with waiting allowed when arriving before $a_i$.

We consider a classic lexicographic objective function: the main goal is to minimize the number of tours $r$, and ties are broken in favor of solutions with a smaller sum of travel times. 
This lexicographic objective function is transformed into a linear objective function by weighting the number of tours $r$ with a large penalty $\lambda = 10n$, and adding all travel times, {\em i.e.}, given a travel time function $\tau$, the cost of a feasible set $S= \{T_1, \ldots, T_r\}$ of $r$ tours is
\begin{equation}
{c}_\tau(S) = \lambda\cdot r + \sum_{T_v\in S} \sum_{(i,j)\in T_v} \tau_{ij}.    
\end{equation}
To avoid numerical instabilities from non-associative operations on floats, we sort the routes of $S$ by their first client to impose a fixed order of evaluation when computing the sum of travel times. 

\section{Scaling and Rounding Schemes}
\label{sec:methodology}

Given a Euclidean instance $I$ and a scaling factor $P=10^k$ with $k \in \mathbb{N}$, we denote by $I^P$ the scaled instance obtained by multiplying all times of $I$ by $P$, and we propose three rounding schemes to handle the conversion of floating-point values to integer values: pessimistic, optimistic, and nearest integer. These rounding schemes are introduced in Sections 3.1 to 3.3. In Section 3.4, we define the scaled-back float cost associated with a feasible solution of $I^P$.

\subsection{Pessimistic Scheme} 

This scheme ensures that any solution feasible for the scaled instance $I^P$ is also feasible for the original instance $I$. This is an \emph{inner approximation} of the feasible region: it makes travel and service times longer while tightening time windows. More precisely, travel and service times are rounded \emph{up} (\texttt{ceil}), {\em i.e.}, 
\[\forall i,j\in C\cup\{d\}, \tau^P_{ij} = \lceil \tau_{ij} \times P \rceil \wedge \forall i\in C, s^P_i = \lceil s_i \times P \rceil\]
and time windows are rounded to the \emph{largest included} integer interval, {\em i.e.}
\[\forall i\in C, a^P_i = \lceil a_i \times P \rceil\wedge b^P_i = \lfloor b_i \times P \rfloor.\]
This tightens constraints to guarantee that feasibility transfers from $I^P$ to $I$. However, some feasible solutions for $I$ may be infeasible for $I^P$. In addition, $I^P$ may be infeasible, {\em e.g.}, when $a^P_i > b^P_i$ for some $i\in C$.

\subsection{Optimistic Scheme} 
This scheme ensures that any solution feasible for the original instance $I$ remains feasible for the scaled instance $I^P$. This is an \emph{outer approximation} of the feasible region: it makes travel and service times shorter while loosening time windows. More precisely, travel and service times are rounded \emph{down} (\texttt{floor}), {\em i.e.},
\[\forall i,j\in C\cup\{d\}, \tau^P_{ij} = \lfloor \tau_{ij} \times P \rfloor \wedge \forall i\in C, s^P_i = \lfloor s_i \times P \rfloor\]
and time windows are rounded to the \emph{smallest enclosing} integer interval, {\em i.e.}, 
\[\forall i\in C, a^P_i = \lfloor a_i \times P \rfloor \wedge b^P_i = \lceil b_i \times P \rceil.\]
This loosens constraints. The optimal cost of $I^P$ (when dividing travel times by $P$) provides a lower bound for the optimal cost of $I$. This scheme, in combination with $P = 10$, is the adopted convention during the 2021 DIMACS VRPTW challenge \cite{DIMACS2021}, and as such has seen wide usage in recent literature.

\subsection{Nearest Integer Scheme}
Instead of controlling feasibility via rounding bias, this scheme aims to minimize average distortion by rounding all values (once multiplied by $P$) to the nearest integer (\texttt{round}).
This scheme offers no strict feasibility guarantees but attempts to minimize numerical distortion on average.

\subsection{Scaled-back Float Cost of Solutions of $I^P$}

Let $S= \{T_1, \ldots, T_r\}$ be a \emph{feasible} solution of a scaled instance $I^P$. To preserve the lexicographic balance between fleet size and travel cost, the vehicle penalty on $I^P$ is set to $\lambda^P = \lambda \cdot P = 10nP$. The integer cost of $S$ on $I^P$ is
\[{c}_{\tau^P}(S) = \lambda^P\cdot r + \sum_{T_v\in S} \sum_{(i,j)\in T_v} \tau^P_{ij}.\]
Dividing by $P$ yields the \emph{scaled-back float cost}:
\[\tilde{c}_{\tau^P}(S) = \lambda\cdot r + \sum_{T_v\in S} \sum_{(i,j)\in T_v} \frac{\tau^P_{ij}}{P}\]
where $\tau^P$ is the travel time function of $I^P$.

When considering the pessimistic scheme, $S$ is always feasible for the original instance $I$. However, when considering the optimistic or nearest schemes, it may happen that a solution feasible for $I^P$ is no longer feasible for the original instance $I$. In this case, the cost of this solution for $I$ is set to infinity, {\em i.e.},
\[
\begin{array}{ll}
c_\tau(S) = + \infty &\mbox{if $S$ violates some time windows for $I$}\\
c_\tau(S) = \lambda\cdot r + \sum_{T_v\in S} \sum_{(i,j)\in T_v} \tau_{ij}
& \mbox{otherwise.}
\end{array}
\]
Whenever $c_\tau(S)<+\infty$, we say that $S$ is \emph{float-feasible}.

\section{Experimental Design}
\label{sec:experiments}

Our experiments address the following research questions:
\begin{enumerate}
    \item How do scaling and rounding affect feasibility transfer from $I^P$ back to $I$?
    \item How much objective distortion is induced by rounding? Does rounding to the nearest integer minimize biases?
    \item How do these effects translate into convergence and final quality for the two considered solvers, HGS and OR-Tools?
\end{enumerate}

In this section, we describe the benchmarks, solvers, experimental setup, and performance measures considered to answer these questions.

\paragraph{Benchmarks.} We consider the 56 Solomon instances with $n=100$, and the Gehring--Homberger instances with $n=600$ and $n=1000$ (60 instances per value of $n$). These instances are grouped in 3 subsets, R (Random), C (Clustered) and RC (Random-Clustered), depending on the topology of their client distribution which has a significant impact on arc cost distribution. These instances belong to the classic SINTEF benchmark with floating point costs and we consider the Best-Known Solutions (BKSs) gathered in \cite{Sintef2008}, denoted $c^*$. For the scaled instances $I^P$, the BKS is the best solution found during all our experiments, and it is denoted $c^{*P}$.

\paragraph{Solvers.} We consider two well-known and open source solvers, {\em i.e.}, PyVRP-HGS \cite{vidalHybridGeneticSearch2022,&PyVRP,PyVRP2024} and OR-Tools \cite{&ortools}. OR-Tools uses a cheapest-insertion start followed by guided local search with a fixed vehicle cost to enforce route-count priority.

Our PyVRP-HGS variant minimizes fleet size during a first stage. Then, it optimizes travel costs with fixed number of vehicles during a second stage. To balance time between stages, Stage~1 uses a custom three-criterion policy: (i) a global Stage~1 time budget (fraction of total run time), (ii) an adaptive per-attempt cap based on observed feasible-attempt durations to bound the final infeasible trial, and (iii) a minimum time reservation for Stage~2.

\paragraph{Experimental setup.} Experiments were executed on a Grid5000 HPC cluster equipped with AMD EPYC 7642 (Zen 2 CPU, \texttt{x86\_64}). The runtime limit ({\it TL}) is set to 120\,s when $n=100$, 900\,s when $n=600$, and 1800\,s when $n=1000$.

We consider the following configurations:
    \begin{itemize}
        \item Scaling factor $P \in \{10^0, 10^1, 10^3\}$,
        \item Rounding scheme in \{Optimistic (Opt), Nearest (Nea), Pessimistic (Pes)\}.
    \end{itemize}
For each configuration, each solver is run 3 times, with 3 different random seeds.

\paragraph{Performance measures.}

To evaluate performance at different runtimes ranging from 0 to the time limit {\it TL}, we consider the evolution of the relative gap to BKS, defined by a function $g:[0, {\it TL}]\rightarrow [0,1]$. Let $S_1, S_2, \ldots, S_k$ be the sequence of incumbent solutions computed at times $t_1, t_2, \ldots, t_k$, respectively, during one run of a solver for a scaled instance $I^P$. Let us first define the relative gap $g^P$ to $c^{*P}$, the BKS of $I^P$:
\begin{itemize}
    \item in order to normalize $g^P$ within $[0,1]$, $g^P(t)$ is set to 1 for each run time $t\in [0,t_1[$ (before a first solution is computed);
    \item then, for each solution $S_i$, the gap is updated, {\em i.e.}, 
for each run time $t\in [t_i,t_{i+1}[$, $g^P(t) = \frac{c_{\tau^P}(S_i)-c^{*P}}{c^{*P}}$.
\end{itemize}
Let us now define the relative gap $g$ to $c^*$, the BKS of the initial instance $I$. In this case, it may happen that $c_\tau(S_i)=+\infty$, because $S_i$ is no longer feasible when considering the cost function $\tau$. Hence, for each solution $S_i$, the gap is updated only if $S_i$ is better than the last computed solution (when evaluated on $\tau$), {\em i.e.}, for each runtime $t\in [t_i,t_{i+1}[$, $g(t) = \min ~\{g(t_{i-1}), \frac{c_\tau(S_i)-c^*}{c^*}\}$ (considering that $t_0=0$).

Finally, to compare relative gap functions, we compute the Area Under the Curve (AUC) for each gap function $g^P$ (resp. $g$). This performance measure is equivalent to the primal-integral measure introduced in \cite{bertholdMeasuringImpactPrimal2013} for MIP solvers.   
More precisely, the normalized AUC corresponds to the area between $g^P(t)$ (resp. $g(t)$) and $y=0$ when the run time $t$ ranges from 0 to the time limit {\it TL}, {\em i.e.}, 
\[
\mathrm{AUC}^P=\frac{1}{\it TL}\int_{0}^{\it TL} g^P(t)\,dt \mbox{ (resp. }\mathrm{AUC}=\frac{1}{\it TL}\int_{0}^{\it TL} g(t)\,dt\mbox{)}.
\]
 When the AUC is equal to 0, we have a perfect run that finds a solution of cost $c^{*P}$ (resp. $c^*$) at time 0; when it is equal to 1, we have a run that is never able to find a solution whose cost is smaller than $2c^{*P}$ (resp. $2c^*$). 

To compare two configurations, we compute the average AUC for each configuration and instance (over the 3 seeds), and then pair average AUCs for each instance. We use paired Wilcoxon tests (for directional differences) and paired TOST-style tests with $\delta=0.01$ (for practical equivalence). All $p$-values are corrected for multiple testing via the Benjamini--Hochberg (BH-FDR) procedure.

\paragraph{Unreported results.} We have also made experiments with a mono-objective variant (where the goal is to minimize the sum of all travel times while ensuring that the number of routes does not exceed some given bound, similar to the DIMACS challenge \cite{DIMACS2021}) and other scaling factors ({\em e.g.}, $P \in \{10^2, 10^5\}$); observed trends were qualitatively similar and are omitted due to the page limit.

\section{Results and Discussion}
\label{sec:results}

\subsection{Feasibility Transfer from $I^P$ to $I$}\label{sec:float-feasibility-transfer}

We first examine whether solutions computed on~$I^P$ remain feasible when re-evaluated on~$I$. During each run, a solver produces a sequence of $k$ incumbents $S_1, \ldots, S_k$. We check every incumbent against the TW constraints of~$I$ and compute the following three metrics: (i)~\emph{first-incumbent feasibility}, {\em i.e.}, the fraction of runs whose initial solution~$S_1$ is float-feasible; (ii)~\emph{all-incumbent feasibility}, {\em i.e.}, for each run the proportion of its $k$ incumbents that are float-feasible, averaged across runs; and (iii)~\emph{last-incumbent feasibility}, {\em i.e.}, the fraction of runs whose final solution~$S_k$ is float-feasible. These metrics respectively capture the initial state, the overall trajectory, and the end-point of the solver's search regarding feasibility transfer. The pessimistic (Pes) scheme guarantees feasibility by construction and yields 100\% for all three metrics; it is therefore omitted.

\begin{figure}[htbp]
  \centering
  \includegraphics[width=\textwidth]{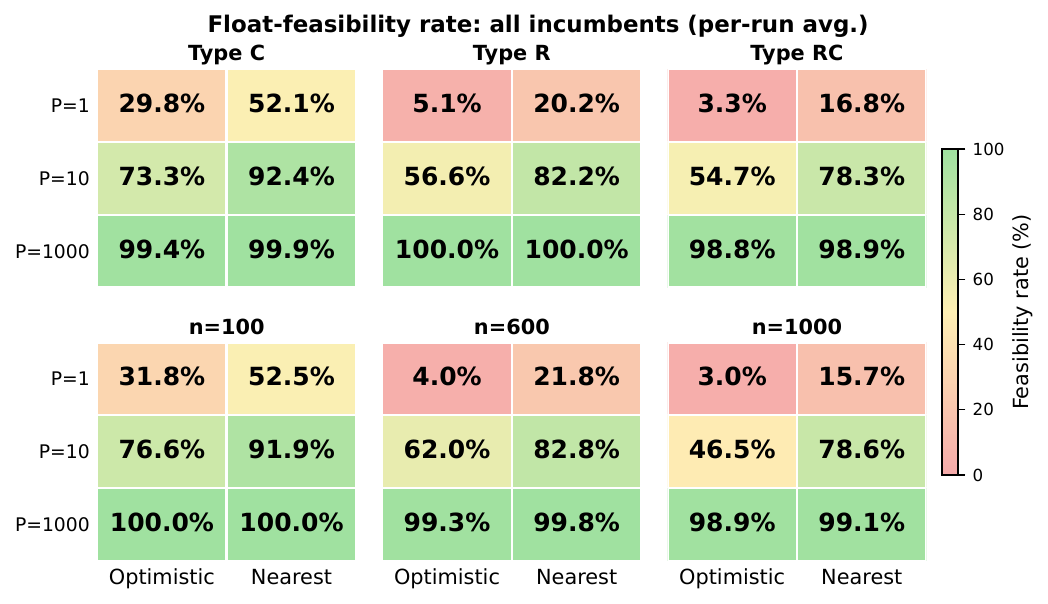}
  \caption{All-incumbent feasibility rate, by rounding scheme and scaling factor, decomposed by instance type (top) and size (bottom). Pessimistic omitted (100\% by construction).}
  \label{fig:feasibility-heatmap}
\end{figure}

Figure~\ref{fig:feasibility-heatmap} displays metric~(ii), decomposed by instance type (R, C, RC) and size~($n$). Several trends emerge. First, all-incumbent feasibility rates improve monotonically with~$P$: at $P{=}1000$ the ratio is larger than 98\% for all subsets, whereas at $P{=}1$ it never exceeds 32\% for Optimistic and 53\% for Nearest. Second, Optimistic is consistently less float-feasible than Nearest, as expected from its more aggressive constraint relaxation. Third, instance topology has a strong impact on float-feasibility when $P\in \{1,10\}$: C-type instances have larger float-feasibility rates than R-type and RC-type instances, for both Optimistic and Nearest. Finally, the number of customers $n$ also has a strong impact on float-feasibility when $P\in \{1,10\}$: the larger $n$, the smaller the float-feasibility rates, {\em e.g.}, for Optimistic with $P{=}1$, float-feasibility decreases from 32\% at $n{=}100$ to 4\% at $n{=}600$ and 3\% at $n{=}1000$.

The first (i) and last (iii) incumbent metrics are not shown for space reasons but follow similar trends, with some nuances that reveal additional dynamics. At low precision, feasibility increases along the search trajectory, {\em e.g.}, for Nearest with $P{=}1$: 22.8\% (first) $\to$ 29.5\% (all) $\to$ 32.2\% (last), indicating that solvers, while not explicitly targeting float-feasibility, progressively drift toward float-feasible regions. This monotonic improvement holds for most type--size combinations, with one notable exception: for C-type instances at $n{=}100$ with Optimistic/$P{=}1$, the first incumbent is more often feasible than the last (80.4\% vs.\ 68.6\%), meaning that cost-driven optimization pushes initially feasible solutions out of the float-feasible region. At $P{=}1000$, all three metrics converge to more than 99\%.

Per-solver analysis reveals a \emph{trajectory crossover} between the two solvers. PyVRP-HGS produces far fewer incumbents per run (${\sim}23$ vs.\ ${\sim}800$ for OR-Tools) but achieves higher first-incumbent feasibility (up to $+22$ percentage points for Optimistic with $P{=}10$). This may be explained by the construction heuristics that generate solutions structurally closer to the float-feasible region. 


These results establish a first key finding: scaling and rounding 
shape float-feasibility throughout the entire search trajectory, with effects modulated by instance topology, size, and solver architecture.

\subsection{Rounding-Induced Cost Estimation Error}

\begin{figure}[htbp]
  \centering
    \includegraphics[width=\textwidth]{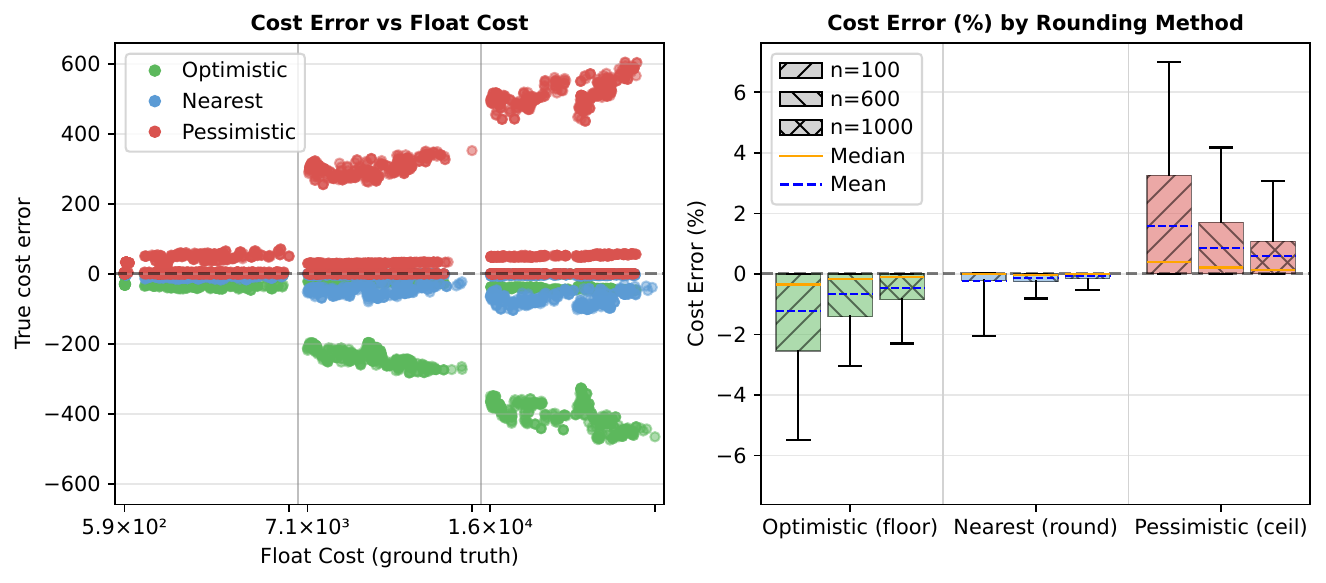}
    \caption{Cost-estimation error induced by integer rounding, shown in absolute and relative forms for the three rounding schemes and scaling factors.}
    \label{fig:cost-comparison-scaling-rounding}
\end{figure}

We now study the rounding-induced cost estimation error by comparing $\tilde{c}_{\tau^P}(S_k)$ (travel times evaluated with $\tau^P$ and rescaled by $1/P$) with $c_\tau(S_k)$ (travel times evaluated with $\tau$), where $S_k$ is the last computed solution. Only runs where $S_k$ is float-feasible are included. Results are presented in Fig.~\ref{fig:cost-comparison-scaling-rounding}. 

In Fig.~\ref{fig:cost-comparison-scaling-rounding} (left), we display a scatter plot of the true cost error $\Delta c(S_k) = \tilde{c}_{\tau^P}(S_k) - c_{\tau}(S_k)$ ($y$-axis) against the ground-truth floating-point cost~$c_{\tau}(S_k)$ ($x$-axis). The $x$-axis uses a piecewise-linear scale giving equal width to each instance-size group; the three resulting segments correspond to $n \in \{100, 600, 1000\}$ (increasing cost ranges from left to right enables this clear separation). Each point represents one run; color indicates the rounding scheme (green = Optimistic, blue = Nearest,  red = Pessimistic). The dashed horizontal line marks zero error; points above the line indicate overestimation of the true Euclidean cost, and points below indicate underestimation. We clearly see several clusters of points corresponding to different rounding schemes and that the absolute error indeed increases with the instance size since larger instances have more edges and thus more rounding operations. Within each size segment, optimistic and pessimistic schemes each exhibit two distinct sub-bands: one near the zero line and one further away. These sub-bands correspond to different scaling factors~$P$: at~$P{=}1$ each rounding operation can displace a distance by up to~${\pm}1$, producing the far-from-zero band; at~$P{=}10$ and~$P{=}1{,}000$ the per-edge error shrinks by a factor of roughly~$10{\times}$ and~$1{,}000{\times}$ respectively, collapsing those runs onto the zero line. Jointly, the scaling factor and instance size explain over~98\% of the variance in absolute error for both the Optimistic and Pessimistic schemes.

In Fig.~\ref{fig:cost-comparison-scaling-rounding} (right), we display grouped boxplots of the relative cost estimation error $\Delta c(S_k) / c_{\tau}(S_k) \times 100$ (in~\%), broken down by rounding scheme ($x$-axis) and instance size (hatch pattern, in increasing order from left to right). Whiskers extend to the full data range (0th--100th percentile); orange lines show the median, blue dashed lines the mean; the $y$-axis is symmetric around zero. On the left panel the absolute error is roughly constant across a given instance size and increases with larger instance sizes, while on the right panel the relative error clearly decreases with increasing instance size: larger instances produce larger total costs, which dilute the fixed absolute rounding error. 

We also observe that the \emph{Nearest} rounding scheme is systematically biased towards underestimation (negative error). We identify two complementary causes and test four hypotheses: (H1)~intrinsic asymmetry of Euclidean-distance fractional parts after scaling, (H2)~optimization-induced selection bias (optimizer's curse), (H3)~possible effect of Python's banker's rounding, and (H4)~attenuation of bias as $P$ increases. Results support H1, H2, and H4, and reject H3.

An edge-level analysis over all 176~instances shows that, at $P{=}1$, nearest rounding rounds down 55.0\% of pairwise distances but rounds up only 41.8\% (the remaining 3.2\% are exact integers), because Euclidean distances from integer coordinates have fractional parts biased towards $[0,\,0.5)$. A one-sample Wilcoxon signed-rank test on per-instance mean rounding errors confirms that this downward bias is highly significant ($W{=}0$, $p<10^{-30}$, $n{=}176$), supporting~H1 independently of any optimization. A cross-rounding re-evaluation experiment then isolates an additional {selection} effect~(H2): when routes optimized under a different scheme (Pessimistic or Optimistic) are merely re-scored with nearest arithmetic, the mean error is about $-0.092\%$; when routes are optimized {under} nearest rounding, the mean error increases to $-0.146\%$, {\em i.e.}, roughly 58\% larger. A paired Wilcoxon test confirms that this gap is significant ($W{=}3\,064$, $p<10^{-170}$, $n{=}1\,056$ pairs at $P{=}1$), consistent with an optimizer's-curse mechanism whereby the solver preferentially selects edges whose rounded costs happen to underestimate the true Euclidean cost. Decomposing the total bias, the intrinsic component accounts for about 63\% and the optimizer's curse for the remaining 37\% at $P{=}1$; this share drops to 14\% at $P{=}10$ and below 2\% at $P{=}1000$, confirming that the bias vanishes with increasing precision~(H4). Finally, an exhaustive comparison over 41.1\,M edge pairs across all instances and scaling factors finds zero edges where banker's rounding and standard half-up disagree~(H3~rejected).

\subsection{Convergence and Performance Under Two Views}

Building on the feasibility and cost-estimation analyses, we now turn to convergence speed and final solution quality. We hypothesize that (H5)~integer-internal performance favors coarser rounding and lower~$P$, because the simplified landscape accelerates search, while (H6)~float-grounded performance favors higher precision and/or pessimistic rounding, because transfer fidelity becomes crucial to approach~$c^*$.

We introduce a \emph{gap-goal reaching ratio} $h_\sigma:[0,{\it TL}]\rightarrow[0,1]$. For a threshold~$\gamma > 0$ and a configuration~$\sigma = (\text{solver},\, n,\, P,\, \text{rounding})$, $h_\sigma(t)$ is the fraction of \emph{all} runs in~$\sigma$ for which~$g(t)\leq\gamma$, where $g$ is the float-grounded gap defined in Section~\ref{sec:experiments}. Every incumbent is checked against the constraints of~$I$; runs that never find a float-feasible solution satisfy $g(t)=1>\gamma$ at all times and count against the ratio. The function $h_\sigma$ is non-decreasing: its asymptotic value $h_\sigma({\it TL})$ simultaneously encodes the transfer-feasibility rate and the solution quality of~$\sigma$, without conditioning on a variable-size subset of runs.

\begin{figure}[htbp]
  \centering
  \includegraphics[width=\textwidth]{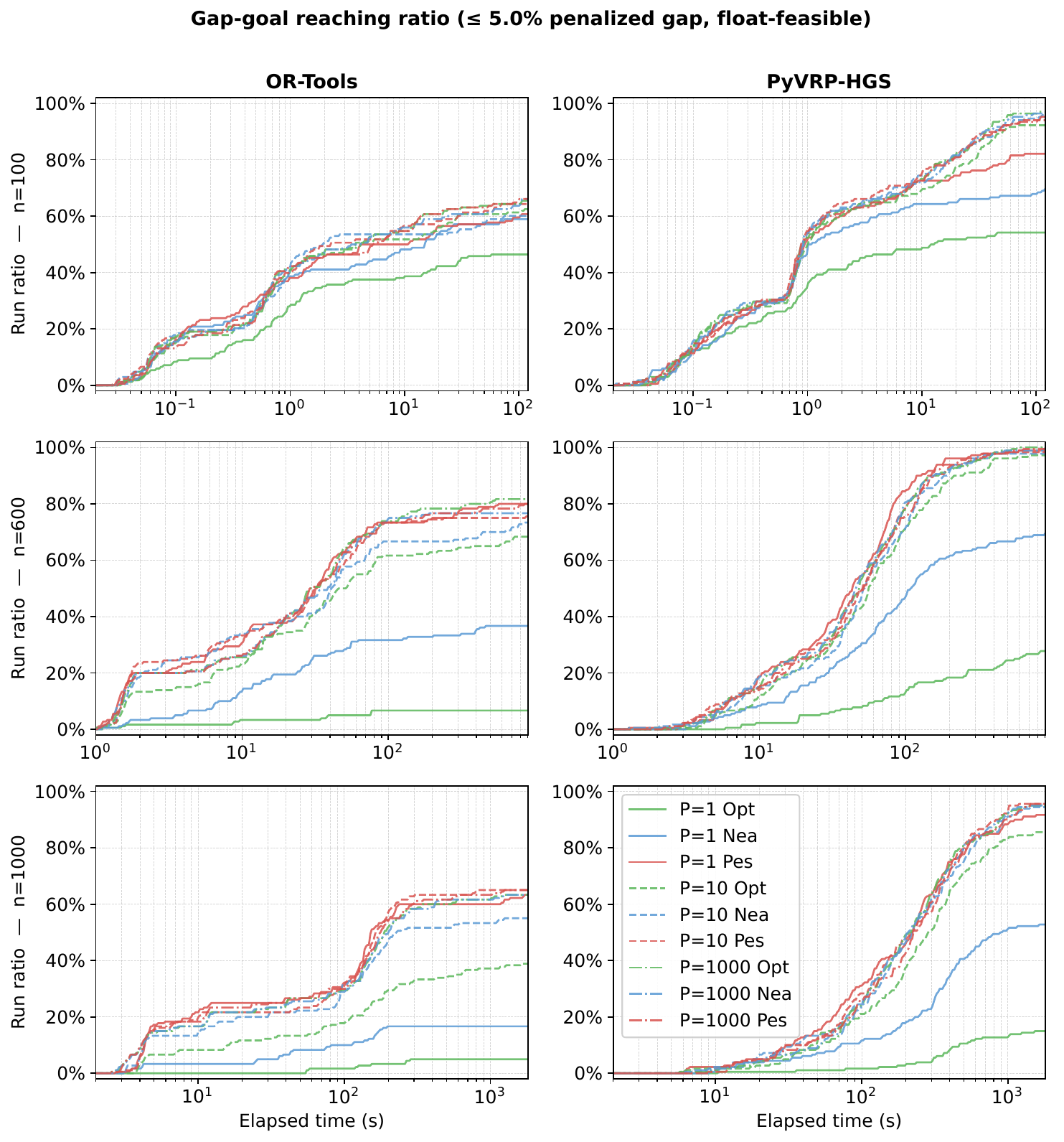}
  \caption{Gap-goal reaching ratio $h_\sigma(t)$ at threshold $\gamma{=}5\%$. Each curve shows the fraction of runs---out of \emph{all} runs in the configuration, including those that never find a float-feasible solution---whose float-grounded gap~$g(t)$ has dropped below~$5\%$ by elapsed time~$t$. Rows: instance sizes; columns: solvers. Colour encodes the rounding scheme (green~= Optimistic, blue~= Nearest, red~= Pessimistic); line style encodes~$P$ ($1$: solid, $10$: dashed, $1000$: dash-dot).}
  \label{fig:gap-goal-5pct-curves}
\end{figure}

Figure~\ref{fig:gap-goal-5pct-curves} plots $h_\sigma(t)$ with $\gamma{=}5\%$ for every configuration. Aggregating on all sizes, OR-Tools reaches only 48--64\%, with curves rising steeply in the first seconds while PyVRP-HGS reaches the goal in 80--87\% of runs, with a more gradual rise. The gap between solvers is a consistent 20--35~percentage points. Within each solver, pessimistic rounding and higher~$P$ yield higher reaching ratios, consistent with Section~\ref{sec:float-feasibility-transfer}, while the worst configuration (Optimistic/$P{=}1$) reaches the lowest ratio across all sizes and solvers. It clearly illustrates that a BKS-centric view of convergence favors high precision and pessimistic rounding.

We complement this analysis with pairwise AUC comparisons, which distill each run's full gap trajectory into a single score. Tables~\ref{tab:pairwise_float_auc_all_combined} and~\ref{tab:pairwise_int_auc_all_combined} present all-size pairwise significance/equivalence matrices, one per metric and solver. Each cell compares row versus column on mean AUC using paired Wilcoxon tests for direction and paired TOST-style tests ($\delta=0.01$) for practical equivalence; all $p$-values are corrected with BH-FDR, yielding adjusted $q$-values. When the directional test is significant ($q_{\text{diff}}<0.05$), the symbol indicates whether the row is better (\texttt{+}) or worse (\texttt{-}); when the equivalence test is significant ($q_{\text{eq}}<0.05$), the symbol is~\texttt{=}. In both cases the number of symbols encodes the strength: one for $q<0.05$, two for $q<0.01$, three for $q<0.001$. Cells marked \texttt{ns} are inconclusive (neither test significant).

\input{tex/tables/2026_02_09_21_03_36_pairwise_matrix_float_auc_all_combined.tex}

\paragraph{Float-grounded AUC (Table~\ref{tab:pairwise_float_auc_all_combined}).} Because this metric assigns a penalty score of $1.0$ to runs lacking a float-feasible final incumbent, configurations with low transfer feasibility are heavily penalized, consistent with the gap-goal reaching ratios of Figure~\ref{fig:gap-goal-5pct-curves}. ($P{=}1$, {Opt}) and ($P{=}1$, {Nea}) are thus significantly worse (\texttt{-{}-{}-}) for both solvers than all higher-precision or pessimistic configurations. In the OR-Tools matrix, ($P{=}1$, {Pes}) already significantly outperforms both (\texttt{+++}), while configurations with $P{\geq}10$ form a tightly equivalent cluster (mostly \texttt{===} or \texttt{ns}). In the PyVRP-HGS matrix the pattern is sharper: ($P{=}1$, {Pes}) is equivalent to the entire $P{=}1000$ block (\texttt{===}), and these jointly dominate everything at $P{\leq}10$ with Opt or Nea. This confirms~H6: float-grounded quality favors either high precision or the pessimistic scheme.

\input{tex/tables/2026_02_09_21_03_36_pairwise_matrix_int_auc_all_combined.tex}

\paragraph{Integer-internal AUC (Table~\ref{tab:pairwise_int_auc_all_combined}).} Because this metric tracks the solver's own objective trajectory on~$I^P$ without any feasibility penalty, it captures a purely solver-centric view of convergence efficiency. The ranking inverts. For both solvers, ($P{=}1$, {Opt}) is the overall winner: it is significantly better (\texttt{+++}) than every other configuration. ($P{=}1$, {Nea}) ranks second, also dominating most $P{\geq}10$ alternatives. Conversely, ($P{=}1$, {Pes}) ranks last, significantly worse (\texttt{-{}-{}-}) than all others; the tightened constraints of the pessimistic scheme slow internal convergence at low precision. For $P{\geq}10$, all rounding schemes become practically equivalent (\texttt{===}), consistent with the idea that the three schemes converge to the same instance as~$P$ grows. This confirms~H5: integer-internal convergence benefits from the coarser landscape of low-precision optimistic rounding.

\paragraph{Synthesis.} To understand why the two tables yield opposite rankings, it is important to recall what each metric fundamentally measures. The float-grounded AUC evaluates convergence toward established BKS on the original floating-point instance~$I$: it uses re-evaluated Euclidean costs, anchors the gap to the SINTEF BKS, and assigns the worst possible score to runs that fail to produce a float-feasible incumbent. This makes it a metric of {benchmark fidelity}, dominated by feasibility transfer. The integer-internal AUC, by contrast, evaluates how efficiently the solver minimizes its own linearized penalized objective $c_{\tau^P}$ on the integer instance~$I^P$: it uses the solver's reported costs directly, computes the reference from observed solver performance, and applies no feasibility penalty. This makes it a metric of {raw convergence efficiency} on the working instance, independent of whether the solution maps back to~$I$.

\input{tex/tables/2026_02_09_21_03_36_auc_gap_feas_summary_all.tex}

Table~\ref{tab:auc_gap_feas_summary_all} crystallizes the resulting tension. The best float-AUC configurations---($P{=}1000$, {Pes}) for OR-Tools and ($P{=}1$, {Pes}) for PyVRP-HGS---are transfer-robust (100\% feasibility), but their float-feasible median gaps to BKS are 11.3\% and 3.8\% respectively. The best integer-AUC configuration is $(P{=}1, \text{Opt})$ for both solvers, with substantially lower gaps among float-feasible runs (median 0.25\% for OR-Tools, 0.03\% for PyVRP-HGS), but transfer coverage is only 8.5\% and 12.7\%. In other words, the few runs that do transfer back under Optimistic/$P{=}1$ happen to be excellent solutions, yet most runs fail the transfer altogether. This ranking inversion is not a contradiction but a direct consequence of the two evaluation perspectives: on the coarse, relaxed landscape of~$I^P$ with low~$P$ and optimistic rounding, the solver converges faster because it explores a simpler search space with a larger feasible region; however, this very simplification means that most final incumbents exploit feasibility margins that do not exist in~$I$, and are therefore invalid when transferred back.

\section{Conclusion}
\label{sec:conclusion}

Scaling and rounding are not neutral implementation details in VRPTW benchmarking; they reshape the feasible region and the fitness landscape, and their impact propagates to  feasibility transfer, cost estimation, and convergence.

Our experiments reveal a clear trade-off. When a solver operates on the integer instance~$I^P$ and fidelity to the original floating-point instance~$I$ is secondary, coarse optimistic rounding at low precision ($P{=}1$, Opt) accelerates integer-internal convergence: it consistently dominates all other configurations in integer AUC for both HGS and OR-Tools. This is coherent with the simplified landscape and enlarged feasible region of the outer approximation. Conversely, when comparison with established BKS on the original instance is the primary goal, higher precision and/or pessimistic rounding is necessary to guarantee transfer feasibility and meaningful float-grounded quality. These two objectives yield a statistically significant ranking inversion between evaluation views.

Beyond this central result, we show that nearest rounding, though intuitive, carries a systematic underestimation bias at low precision, driven by the distribution of Euclidean-distance fractional parts (63\% of the effect) and amplified by an optimizer's-curse mechanism (37\%), both vanishing as~$P$ grows.

In practice, we offer three recommendations: (i)~explicitly state the scaling factor, rounding scheme, and evaluation space when reporting VRPTW results; (ii)~prefer the pessimistic scheme or $P \geq 10^3$ when reproducible comparison to float-based BKS is important; (iii)~when integer-space performance is the sole criterion, low-precision optimistic rounding is a valid and efficient choice. This work is a first step toward updated VRPTW benchmark practices, tooling and instance definitions that better support reproducibility and fair comparison between solvers.

\begin{credits}
\subsubsection{\ackname}
This work was funded as part of the \href{https://anr.fr/Project-ANR-22-CE22-0016}{ANR project MAMUT (ANR-22-CE22-0016)}, France. Special thanks to Romain Fontaine for his valuable insights. Experiments presented in this paper were carried out using the \href{https://www.grid5000.fr}{Grid'5000 testbed}, supported by a scientific interest group hosted by \href{https://www.inria.fr/en}{Inria} and including \href{https://www.cnrs.fr/en}{CNRS}, \href{https://www.renater.fr/en/accueil-english/}{RENATER} and several Universities as well as other organizations\footnote{See \href{https://www.grid5000.fr}{https://www.grid5000.fr}}. 

\subsubsection{\discintname}
The authors have no competing interests to declare.
\end{credits}

%
%
\bibliographystyle{splncs04}
\bibliography{biblio}
\end{document}

%% file: tex/tables/2026_02_09_21_03_36_pairwise_matrix_float_auc_all_combined.tex
\begin{table}[htbp]
\centering
\normalsize
\caption{Per solver pairwise Float-grounded normalized AUC outcomes}
\label{tab:pairwise_float_auc_all_combined}
\textbf{OR-Tools}\\
\begin{tabular*}{\linewidth}{@{\extracolsep{\fill}}lccccccccc@{}}
\toprule
 & 1-Opt & 1-Nea & 1-Pes & 10-Opt & 10-Nea & 10-Pes & 1000-Opt & 1000-Nea & 1000-Pes \\
\midrule
1-Opt &  & \texttt{ns} & \cellcolor{matworse3}\texttt{---} & \cellcolor{matworse3}\texttt{---} & \cellcolor{matworse3}\texttt{---} & \cellcolor{matworse3}\texttt{---} & \cellcolor{matworse3}\texttt{---} & \cellcolor{matworse3}\texttt{---} & \cellcolor{matworse3}\texttt{---} \\
1-Nea & \texttt{ns} &  & \cellcolor{matworse3}\texttt{---} & \texttt{ns} & \cellcolor{matworse3}\texttt{---} & \cellcolor{matworse3}\texttt{---} & \cellcolor{matworse3}\texttt{---} & \cellcolor{matworse3}\texttt{---} & \cellcolor{matworse3}\texttt{---} \\
1-Pes & \cellcolor{matbetter3}\texttt{+++} & \cellcolor{matbetter3}\texttt{+++} &  & \cellcolor{matbetter3}\texttt{+++} & \cellcolor{matworse3}\texttt{---} & \cellcolor{matworse3}\texttt{---} & \cellcolor{matworse3}\texttt{---} & \cellcolor{matworse3}\texttt{---} & \cellcolor{matworse3}\texttt{---} \\
10-Opt & \cellcolor{matbetter3}\texttt{+++} & \texttt{ns} & \cellcolor{matworse3}\texttt{---} &  & \cellcolor{matworse3}\texttt{---} & \cellcolor{matworse3}\texttt{---} & \cellcolor{matworse3}\texttt{---} & \cellcolor{matworse3}\texttt{---} & \cellcolor{matworse3}\texttt{---} \\
10-Nea & \cellcolor{matbetter3}\texttt{+++} & \cellcolor{matbetter3}\texttt{+++} & \cellcolor{matbetter3}\texttt{+++} & \cellcolor{matbetter3}\texttt{+++} &  & \texttt{ns} & \cellcolor{matworse2}\texttt{--} & \texttt{ns} & \texttt{ns} \\
10-Pes & \cellcolor{matbetter3}\texttt{+++} & \cellcolor{matbetter3}\texttt{+++} & \cellcolor{matbetter3}\texttt{+++} & \cellcolor{matbetter3}\texttt{+++} & \texttt{ns} &  & \texttt{ns} & \texttt{ns} & \cellcolor{matequiv3}\texttt{===} \\
1000-Opt & \cellcolor{matbetter3}\texttt{+++} & \cellcolor{matbetter3}\texttt{+++} & \cellcolor{matbetter3}\texttt{+++} & \cellcolor{matbetter3}\texttt{+++} & \cellcolor{matbetter2}\texttt{++} & \texttt{ns} &  & \cellcolor{matequiv3}\texttt{===} & \texttt{ns} \\
1000-Nea & \cellcolor{matbetter3}\texttt{+++} & \cellcolor{matbetter3}\texttt{+++} & \cellcolor{matbetter3}\texttt{+++} & \cellcolor{matbetter3}\texttt{+++} & \texttt{ns} & \texttt{ns} & \cellcolor{matequiv3}\texttt{===} &  & \texttt{ns} \\
1000-Pes & \cellcolor{matbetter3}\texttt{+++} & \cellcolor{matbetter3}\texttt{+++} & \cellcolor{matbetter3}\texttt{+++} & \cellcolor{matbetter3}\texttt{+++} & \texttt{ns} & \cellcolor{matequiv3}\texttt{===} & \texttt{ns} & \texttt{ns} &  \\
\bottomrule
\end{tabular*}
\vspace{0.5em}
\textbf{PyVRP-HGS}\\
\begin{tabular*}{\linewidth}{@{\extracolsep{\fill}}lccccccccc@{}}
\toprule
 & 1-Opt & 1-Nea & 1-Pes & 10-Opt & 10-Nea & 10-Pes & 1000-Opt & 1000-Nea & 1000-Pes \\
\midrule
1-Opt &  & \texttt{ns} & \cellcolor{matworse3}\texttt{---} & \cellcolor{matworse3}\texttt{---} & \cellcolor{matworse3}\texttt{---} & \cellcolor{matworse3}\texttt{---} & \cellcolor{matworse3}\texttt{---} & \cellcolor{matworse3}\texttt{---} & \cellcolor{matworse3}\texttt{---} \\
1-Nea & \texttt{ns} &  & \cellcolor{matworse3}\texttt{---} & \cellcolor{matworse3}\texttt{---} & \cellcolor{matworse3}\texttt{---} & \cellcolor{matworse3}\texttt{---} & \cellcolor{matworse3}\texttt{---} & \cellcolor{matworse3}\texttt{---} & \cellcolor{matworse3}\texttt{---} \\
1-Pes & \cellcolor{matbetter3}\texttt{+++} & \cellcolor{matbetter3}\texttt{+++} &  & \cellcolor{matbetter3}\texttt{+++} & \cellcolor{matworse3}\texttt{---} & \cellcolor{matequiv3}\texttt{===} & \cellcolor{matequiv3}\texttt{===} & \cellcolor{matequiv2}\texttt{==} & \cellcolor{matequiv3}\texttt{===} \\
10-Opt & \cellcolor{matbetter3}\texttt{+++} & \cellcolor{matbetter3}\texttt{+++} & \cellcolor{matworse3}\texttt{---} &  & \cellcolor{matworse3}\texttt{---} & \cellcolor{matworse3}\texttt{---} & \cellcolor{matworse3}\texttt{---} & \cellcolor{matworse3}\texttt{---} & \cellcolor{matworse3}\texttt{---} \\
10-Nea & \cellcolor{matbetter3}\texttt{+++} & \cellcolor{matbetter3}\texttt{+++} & \cellcolor{matbetter3}\texttt{+++} & \cellcolor{matbetter3}\texttt{+++} &  & \cellcolor{matworse3}\texttt{---} & \cellcolor{matworse3}\texttt{---} & \cellcolor{matworse3}\texttt{---} & \cellcolor{matworse3}\texttt{---} \\
10-Pes & \cellcolor{matbetter3}\texttt{+++} & \cellcolor{matbetter3}\texttt{+++} & \cellcolor{matequiv3}\texttt{===} & \cellcolor{matbetter3}\texttt{+++} & \cellcolor{matbetter3}\texttt{+++} &  & \cellcolor{matequiv3}\texttt{===} & \cellcolor{matequiv2}\texttt{==} & \cellcolor{matequiv3}\texttt{===} \\
1000-Opt & \cellcolor{matbetter3}\texttt{+++} & \cellcolor{matbetter3}\texttt{+++} & \cellcolor{matequiv3}\texttt{===} & \cellcolor{matbetter3}\texttt{+++} & \cellcolor{matbetter3}\texttt{+++} & \cellcolor{matequiv3}\texttt{===} &  & \cellcolor{matequiv3}\texttt{===} & \cellcolor{matequiv3}\texttt{===} \\
1000-Nea & \cellcolor{matbetter3}\texttt{+++} & \cellcolor{matbetter3}\texttt{+++} & \cellcolor{matequiv2}\texttt{==} & \cellcolor{matbetter3}\texttt{+++} & \cellcolor{matbetter3}\texttt{+++} & \cellcolor{matequiv2}\texttt{==} & \cellcolor{matequiv3}\texttt{===} &  & \cellcolor{matequiv3}\texttt{===} \\
1000-Pes & \cellcolor{matbetter3}\texttt{+++} & \cellcolor{matbetter3}\texttt{+++} & \cellcolor{matequiv3}\texttt{===} & \cellcolor{matbetter3}\texttt{+++} & \cellcolor{matbetter3}\texttt{+++} & \cellcolor{matequiv3}\texttt{===} & \cellcolor{matequiv3}\texttt{===} & \cellcolor{matequiv3}\texttt{===} &  \\
\bottomrule
\end{tabular*}
\vspace{0.35em}
\begin{minipage}{\linewidth}
\footnotesize
\textit{Notes:} Cell compares row vs column. Lower AUC is better. Symbols: \texttt{+}/\texttt{++}/\texttt{+++}: significantly better, \texttt{-}/\texttt{-{}-}/\texttt{-{}-{}-}: significantly worse, \texttt{=}/\texttt{==}/\texttt{===}: significant practical equivalence, \texttt{ns}: inconclusive. 
\end{minipage}
\end{table}

%% file: tex/tables/2026_02_09_21_03_36_pairwise_matrix_int_auc_all_combined.tex
\begin{table}[htbp]
\centering
\normalsize
\caption{Per solver pairwise Integer-internal normalized AUC outcomes}
\label{tab:pairwise_int_auc_all_combined}
\textbf{OR-Tools}\\
\begin{tabular*}{\linewidth}{@{\extracolsep{\fill}}lccccccccc@{}}
\toprule
 & 1-Opt & 1-Nea & 1-Pes & 10-Opt & 10-Nea & 10-Pes & 1000-Opt & 1000-Nea & 1000-Pes \\
\midrule
1-Opt &  & \cellcolor{matbetter3}\texttt{+++} & \cellcolor{matbetter3}\texttt{+++} & \cellcolor{matbetter3}\texttt{+++} & \cellcolor{matbetter3}\texttt{+++} & \cellcolor{matbetter3}\texttt{+++} & \cellcolor{matbetter3}\texttt{+++} & \cellcolor{matbetter3}\texttt{+++} & \cellcolor{matbetter3}\texttt{+++} \\
1-Nea & \cellcolor{matworse3}\texttt{---} &  & \cellcolor{matbetter3}\texttt{+++} & \cellcolor{matbetter2}\texttt{++} & \cellcolor{matbetter1}\texttt{+} & \cellcolor{matbetter3}\texttt{+++} & \cellcolor{matbetter3}\texttt{+++} & \cellcolor{matbetter3}\texttt{+++} & \cellcolor{matbetter3}\texttt{+++} \\
1-Pes & \cellcolor{matworse3}\texttt{---} & \cellcolor{matworse3}\texttt{---} &  & \cellcolor{matworse3}\texttt{---} & \cellcolor{matworse3}\texttt{---} & \cellcolor{matworse3}\texttt{---} & \cellcolor{matworse3}\texttt{---} & \cellcolor{matworse3}\texttt{---} & \cellcolor{matworse3}\texttt{---} \\
10-Opt & \cellcolor{matworse3}\texttt{---} & \cellcolor{matworse2}\texttt{--} & \cellcolor{matbetter3}\texttt{+++} &  & \cellcolor{matequiv3}\texttt{===} & \cellcolor{matbetter3}\texttt{+++} & \cellcolor{matbetter2}\texttt{++} & \cellcolor{matbetter3}\texttt{+++} & \cellcolor{matbetter3}\texttt{+++} \\
10-Nea & \cellcolor{matworse3}\texttt{---} & \cellcolor{matworse1}\texttt{-} & \cellcolor{matbetter3}\texttt{+++} & \cellcolor{matequiv3}\texttt{===} &  & \cellcolor{matbetter3}\texttt{+++} & \cellcolor{matequiv3}\texttt{===} & \cellcolor{matbetter1}\texttt{+} & \cellcolor{matbetter1}\texttt{+} \\
10-Pes & \cellcolor{matworse3}\texttt{---} & \cellcolor{matworse3}\texttt{---} & \cellcolor{matbetter3}\texttt{+++} & \cellcolor{matworse3}\texttt{---} & \cellcolor{matworse3}\texttt{---} &  & \cellcolor{matequiv3}\texttt{===} & \cellcolor{matequiv3}\texttt{===} & \cellcolor{matequiv3}\texttt{===} \\
1000-Opt & \cellcolor{matworse3}\texttt{---} & \cellcolor{matworse3}\texttt{---} & \cellcolor{matbetter3}\texttt{+++} & \cellcolor{matworse2}\texttt{--} & \cellcolor{matequiv3}\texttt{===} & \cellcolor{matequiv3}\texttt{===} &  & \cellcolor{matequiv3}\texttt{===} & \cellcolor{matequiv3}\texttt{===} \\
1000-Nea & \cellcolor{matworse3}\texttt{---} & \cellcolor{matworse3}\texttt{---} & \cellcolor{matbetter3}\texttt{+++} & \cellcolor{matworse3}\texttt{---} & \cellcolor{matworse1}\texttt{-} & \cellcolor{matequiv3}\texttt{===} & \cellcolor{matequiv3}\texttt{===} &  & \cellcolor{matequiv3}\texttt{===} \\
1000-Pes & \cellcolor{matworse3}\texttt{---} & \cellcolor{matworse3}\texttt{---} & \cellcolor{matbetter3}\texttt{+++} & \cellcolor{matworse3}\texttt{---} & \cellcolor{matworse1}\texttt{-} & \cellcolor{matequiv3}\texttt{===} & \cellcolor{matequiv3}\texttt{===} & \cellcolor{matequiv3}\texttt{===} &  \\
\bottomrule
\end{tabular*}
\vspace{0.5em}
\textbf{PyVRP-HGS}\\
\begin{tabular*}{\linewidth}{@{\extracolsep{\fill}}lccccccccc@{}}
\toprule
 & 1-Opt & 1-Nea & 1-Pes & 10-Opt & 10-Nea & 10-Pes & 1000-Opt & 1000-Nea & 1000-Pes \\
\midrule
1-Opt &  & \cellcolor{matbetter3}\texttt{+++} & \cellcolor{matbetter3}\texttt{+++} & \cellcolor{matbetter3}\texttt{+++} & \cellcolor{matbetter3}\texttt{+++} & \cellcolor{matbetter3}\texttt{+++} & \cellcolor{matbetter3}\texttt{+++} & \cellcolor{matbetter3}\texttt{+++} & \cellcolor{matbetter3}\texttt{+++} \\
1-Nea & \cellcolor{matworse3}\texttt{---} &  & \cellcolor{matbetter3}\texttt{+++} & \cellcolor{matbetter3}\texttt{+++} & \cellcolor{matbetter3}\texttt{+++} & \cellcolor{matbetter3}\texttt{+++} & \cellcolor{matbetter3}\texttt{+++} & \cellcolor{matbetter3}\texttt{+++} & \cellcolor{matbetter3}\texttt{+++} \\
1-Pes & \cellcolor{matworse3}\texttt{---} & \cellcolor{matworse3}\texttt{---} &  & \cellcolor{matworse3}\texttt{---} & \cellcolor{matworse3}\texttt{---} & \cellcolor{matworse3}\texttt{---} & \cellcolor{matworse3}\texttt{---} & \cellcolor{matworse3}\texttt{---} & \cellcolor{matworse3}\texttt{---} \\
10-Opt & \cellcolor{matworse3}\texttt{---} & \cellcolor{matworse3}\texttt{---} & \cellcolor{matbetter3}\texttt{+++} &  & \cellcolor{matbetter3}\texttt{+++} & \cellcolor{matbetter3}\texttt{+++} & \cellcolor{matbetter1}\texttt{+} & \cellcolor{matbetter1}\texttt{+} & \cellcolor{matbetter3}\texttt{+++} \\
10-Nea & \cellcolor{matworse3}\texttt{---} & \cellcolor{matworse3}\texttt{---} & \cellcolor{matbetter3}\texttt{+++} & \cellcolor{matworse3}\texttt{---} &  & \cellcolor{matbetter1}\texttt{+} & \cellcolor{matequiv3}\texttt{===} & \cellcolor{matequiv3}\texttt{===} & \cellcolor{matequiv3}\texttt{===} \\
10-Pes & \cellcolor{matworse3}\texttt{---} & \cellcolor{matworse3}\texttt{---} & \cellcolor{matbetter3}\texttt{+++} & \cellcolor{matworse3}\texttt{---} & \cellcolor{matworse1}\texttt{-} &  & \cellcolor{matworse1}\texttt{-} & \cellcolor{matequiv3}\texttt{===} & \cellcolor{matequiv3}\texttt{===} \\
1000-Opt & \cellcolor{matworse3}\texttt{---} & \cellcolor{matworse3}\texttt{---} & \cellcolor{matbetter3}\texttt{+++} & \cellcolor{matworse1}\texttt{-} & \cellcolor{matequiv3}\texttt{===} & \cellcolor{matbetter1}\texttt{+} &  & \cellcolor{matequiv3}\texttt{===} & \cellcolor{matequiv3}\texttt{===} \\
1000-Nea & \cellcolor{matworse3}\texttt{---} & \cellcolor{matworse3}\texttt{---} & \cellcolor{matbetter3}\texttt{+++} & \cellcolor{matworse1}\texttt{-} & \cellcolor{matequiv3}\texttt{===} & \cellcolor{matequiv3}\texttt{===} & \cellcolor{matequiv3}\texttt{===} &  & \cellcolor{matequiv3}\texttt{===} \\
1000-Pes & \cellcolor{matworse3}\texttt{---} & \cellcolor{matworse3}\texttt{---} & \cellcolor{matbetter3}\texttt{+++} & \cellcolor{matworse3}\texttt{---} & \cellcolor{matequiv3}\texttt{===} & \cellcolor{matequiv3}\texttt{===} & \cellcolor{matequiv3}\texttt{===} & \cellcolor{matequiv3}\texttt{===} &  \\
\bottomrule
\end{tabular*}
\vspace{0.35em}
\begin{minipage}{\linewidth}
\footnotesize
\textit{Notes:} Similar semantics as Table~\ref{tab:pairwise_float_auc_all_combined}.
\end{minipage}
\end{table}

%% file: tex/tables/2026_02_09_21_03_36_auc_gap_feas_summary_all.tex
\begin{table}[htbp]
\centering
\normalsize
\caption{All-size summary: best configuration by AUC with final-gap quality}
\label{tab:auc_gap_feas_summary_all}
\begin{tabular}{l l l r r r}
\toprule
Solver & Metric & Best config & Mean gap (\%) & Median gap (\%) & Feas. (\%) \\
\midrule
OR-Tools & float AUC & 1000-Pes & 9.96 & 11.33 & 100.00 \\
OR-Tools & integer AUC & 1-Opt & 1.00 & 0.25 & 8.52 \\
PyVRP-HGS & float AUC & 1-Pes & 5.50 & 3.80 & 100.00 \\
PyVRP-HGS & integer AUC & 1-Opt & 0.82 & 0.03 & 12.69 \\
\bottomrule
\end{tabular}
\vspace{0.35em}
\begin{minipage}{\linewidth}
\footnotesize
\textit{Notes:} Best configuration is selected by lowest mean AUC on all runs for the given solver and metric. Gap statistics are computed on float-feasible runs only; Feas. is the float-feasible percentage over all runs of that best configuration.
\end{minipage}
\end{table}